\documentclass[10pt]{amsart}
\usepackage[cp1251]{inputenc}
\usepackage[english]{babel}
\usepackage{amsmath}
\usepackage{amssymb}
\usepackage{amsfonts}

\usepackage[linktocpage=true, colorlinks=true, linkcolor=blue, citecolor=blue, urlcolor=blue]{hyperref}

\begin{document}
\renewcommand{\refname}{References}

\thispagestyle{empty}

\title[New Simple Method of Expansion of Iterated Ito
Stochastic Integrals]
{New Simple Method of Expansion of Iterated
Ito Stochastic integrals of Multiplicity 2 Based on Expansion of 
the Brownian Motion Using Legendre Polynomials and Trigonometric 
Functions}
\author[D.F. Kuznetsov]{Dmitriy F. Kuznetsov}
\address{Dmitriy Feliksovich Kuznetsov
\newline\hphantom{iii} Peter the Great Saint-Petersburg Polytechnic University,
\newline\hphantom{iii} Polytechnicheskaya ul., 29,
\newline\hphantom{iii} 195251, Saint-Petersburg, Russia}%
\email{sde\_kuznetsov@inbox.ru}

\maketitle {\small
\begin{quote}
\noindent{\sc Abstract. } 
The article is devoted to the 
expansion of iterated Ito
stochastic integrals of second multiplicity based on expansion of the
Brownian motion (standard Wiener process) using 
complete orthonormal systems of functions in the space
$L_2([t, T]).$ The cases of Le\-gen\-dre polynomials 
and trigonometric functions are considered in details.
We obtained a new re\-pre\-sen\-ta\-ti\-on of the Levy stochastic area based on 
the Legendre polynomials. This re\-pre\-sen\-ta\-ti\-on was first
derived in the author's work \cite{300a} (1997).
In this article, we obtain the mentioned representation 
by a simpler method compared to \cite{300a} (1997).
Also, we get the polynomial re\-pre\-sen\-ta\-ti\-on of the Levy
stochastic area using the method
of expansion of iterated Ito stochastic
integrals based on generalized multple Fourier
series. The polynomial re\-pre\-sen\-ta\-ti\-on of the Levy stochastic area
has more simple form in comparison with the classical trigonometric
re\-pre\-sen\-ta\-ti\-on of the Levy stochastic area.
The convergence in the mean of degree $2n$ $(n\in\mathbb{N})$
as well as the convergence 
with probability 1 for approximations of the Levy stochastic area are proved.
The results of the article can be applied to the numerical
solution of Ito stochastic differential equations
as well as to the numerical approximation of mild
solution for non-commutative semi\-li\-ne\-ar stochastic partial differential 
equations.

\medskip
\end{quote}
}

\vspace{12mm}

\setlength{\baselineskip}{1.5em}

\tableofcontents

\setlength{\baselineskip}{1.2em}

\vspace{5mm}

\section{Introduction}

\vspace{5mm}

Let $(\Omega,$ ${\rm F},$ ${\sf P})$ be a complete probability space, let 
$\{{\rm F}_t, t\in[0,T]\}$ be a nondecreasing right-continous 
family of $\sigma$-algebras of ${\rm F},$
and let ${\bf w}_t$ be a standard $m$-dimensional Wiener 
stochastic process, which is
${\rm F}_t$-measurable for any $t\in[0, T].$ We assume that the components
${\bf w}_{t}^{(i)}$ $(i=1,\ldots,m)$ of this process are independent. 
Consider
an Ito stochastic differential equation in the integral form  

\begin{equation}
\label{1.5.2}
{\bf x}_t={\bf x}_0+\int\limits_0^t {\bf a}({\bf x}_{\tau},\tau)d\tau+
\int\limits_0^t B({\bf x}_{\tau},\tau)d{\bf w}_{\tau},\ \ \
{\bf x}_0={\bf x}(0,\omega),\ \ \ \omega\in\Omega.
\end{equation}

\vspace{3mm}
\noindent
Here ${\bf x}_t$ is some $n$-dimensional stochastic process 
satisfying the equation (\ref{1.5.2}). 
The nonrandom functions ${\bf a}: \mathbb{R}^n\times[0, T]\to\mathbb{R}^n$,
$B: \mathbb{R}^n\times[0, T]\to \mathbb{R}^{n\times m}$
guarantee the existence and uniqueness up to stochastic 
equivalence of a solution
of the equation (\ref{1.5.2}) \cite{1}. The second integral 
on the right-hand side of (\ref{1.5.2}) is 
interpreted as an Ito stochastic integral.
Let ${\bf x}_0$ be an $n$-dimensional random variable, which is 
${\rm F}_0$-measurable and 
${\sf M}\{\left|{\bf x}_0\right|^2\}<\infty$ 
(${\sf M}$ denotes a mathematical expectation).
We assume that
${\bf x}_0$ and ${\bf w}_t-{\bf w}_0$ are independent when $t>0.$

One of the effective approaches 
to the numerical integration of 
Ito stochastic differential equations
is an approach based on the Taylor--Ito 
expansion \cite{2}-\cite{3}.
The most important feature of the Taylor--Ito
expansion is a presence in this expansion of the so-called iterated 
Ito stochastic integrals, which play the key 
role for solving the 
problem of numerical integration of Ito stochastic differential equations
and have the 
following form

\begin{equation}
\label{ito}
J[\psi^{(k)}]_{T,t}=\int\limits_t^T\psi_k(t_k) \ldots \int\limits_t^{t_{2}}
\psi_1(t_1) d{\bf w}_{t_1}^{(i_1)}\ldots
d{\bf w}_{t_k}^{(i_k)}\ \ \ (i_1,\ldots,i_k = 0, 1,\ldots,m),
\end{equation}

\vspace{3mm}
\noindent
where $\psi_1(\tau),\ldots,\psi_k(\tau)$ are nonrandom
functions 
on $[t,T],$ ${\bf w}_{\tau}^{(i)}$
($i=1,\ldots,m$) are independent standard Wiener processes, and
${\bf w}_{\tau}^{(0)}=\tau$.

In this article, we pay a special 
attention to the case $k=2,$\ \ $i_1, i_2=1,\ldots,m$,\ \
$\psi_1(\tau),$ $\psi_2(\tau)\equiv 1$. 
This case corresponds to 
the so-called Milstein method \cite{Mi2}, \cite{3} for
the numerical integration of Ito  stochastic differential equations. 
It is well known that the Milstein method 
has the order 1.0 of strong convergence under 
the specific conditions \cite{Mi2}, \cite{3}.

The Milstein method has the following form \cite{Mi2}, \cite{3}

\vspace{1mm}
$$
{\bf y}_{p+1}={\bf y}_p+\sum_{i_{1}=1}^{m}B_{i_{1}}
I_{\tau_{p+1},\tau_p}^{(i_{1})}+\Delta{\bf a}
+\sum_{i_{1},i_{2}=1}^{m}G_{i_{1}}
B_{i_{2}}\hat I_{\tau_{p+1},\tau_p}^{(i_{1}i_{2})},
$$

\vspace{4mm}
\noindent
where $\Delta=T/N$ $(N>1)$ is a constant (for simplicity)
step of integration,\
$\tau_p=p\Delta$ $(p=0, 1,\ldots,N)$,\

$$
G_i = \sum^ {n} _ {j=1} B_{ji} ({\bf x}, t)
{\partial  \over \partial {\bf x}_j}\ \ \ 
(i=1,\ldots,m),
$$

\vspace{3mm}
\noindent
$B_i$ is the $i$th column of the matrix function $B$ and $B_{ij}$  
is the $ij$th element
of the matrix function 
$B$,  ${\bf a}_i$ is the $i$th element of the 
vector function ${\bf a}$, and ${\bf x}_i$ 
is the
$i$th element of the column ${\bf x}$, 
the columns 
$B_{i_{1}},$\ ${\bf a},$\ $G_{i_{1}}B_{i_{2}}$
are calculated in the point $({\bf y}_p,p),$

$$
I_{\tau_{p+1},\tau_p}^{(i_{1})}=
\int\limits_{\tau_p}^{\tau_{p+1}}d{\bf w}_{\tau}^{(i_1)},
$$

\vspace{3mm}
\noindent
$\hat I_{\tau_{p+1},\tau_p}^{(i_{1}i_{2})}$ is an approximation of the
following iterated Ito stochastic integral

$$
I_{\tau_{p+1},\tau_p}^{(i_{1}i_{2})}=
\int\limits_{\tau_p}^{\tau_{p+1}}\int\limits_{\tau_p}^{s}
d{\bf w}_{\tau}^{(i_1)}d{\bf w}_{s}^{(i_2)}.
$$

\vspace{3mm}

The Levy stochastic area $A^{(i_1i_2)}_{T,t}$ is defined as follows
\cite{1c}

\vspace{2mm}
$$
A_{T,t}^{(i_1 i_2)}=
\frac{1}{2}\left(I_{T,t}^{(i_1 i_2)}-I_{T,t}^{(i_2 i_1)}\right).
$$

\vspace{4mm}

It is clear that

\begin{equation}
\label{2121}
I_{T,t}^{(i_1 i_2)}=\frac{1}{2}I_{T,t}^{(i_1)}I_{T,t}^{(i_2)}+
A_{T,t}^{(i_1 i_2)}\ \ \ \hbox{w. p. 1},
\end{equation}

\vspace{3mm}
\noindent
where w.~p.~1 means with probability 1, $i_1\ne i_2.$

The relation (\ref{2121}) implies that the problem of numerical
simulation of the iterated Ito stochastic integral 
$I_{T,t}^{(i_1 i_2)}$ is equivalent to the problem of numerical
simulation of the Levy stochastic area.

There are some methods for representation of the Levy stochastic area
(see, for example, \cite{2}--\cite{3}). In this article, we consider a new
representation of the Levy stochastic area based
on the Legendre polynomials. This representation is simpler than its
existing analogue based on the Karhunen--Loeve expansion of the 
Brownian bridge process \cite{Mi2} (also see \cite{2}).

\vspace{5mm}

\section{Method of Expansion of Iterated Ito stochastic integrals
Based on Generalized Multiple
Fourier Series}

\vspace{5mm}

Consider the
iterated Ito stochastic integrals (\ref{ito}) and
define the following function on the hypercube $[t, T]^k$

\vspace{-1mm}
$$
K(t_1,\ldots,t_k)=
\begin{cases}
\psi_1(t_1)\ldots \psi_k(t_k),\ &\ t_1<\ldots<t_k\\
~\\
~\\
0,\ &\ \hbox{\rm otherwise}
\end{cases},\ \ \ t_1,\ldots,t_k\in[t, T],\ \ \ k\ge 2,
$$

\vspace{4mm}
\noindent
and 
$K(t_1)\equiv\psi_1(t_1)$ for $t_1\in[t, T].$ 
Here we suppose that $\psi_1(\tau),\ldots,\psi_k(\tau)\in L_2([t, T]).$

Assume that $\{\phi_j(x)\}_{j=0}^{\infty}$
is a complete orthonormal system of functions in 
the space $L_2([t, T])$. 
The function $K(t_1,\ldots,t_k)$ belongs to the space $L_2([t, T]^k)$. 
At this situation it is well known that the generalized 
multiple Fourier series 
of $K(t_1,\ldots,t_k)\in L_2([t, T]^k)$ is converging 
to $K(t_1,\ldots,t_k)$ in the hypercube $[t, T]^k$ in 
the mean-square sense, i.e.

\begin{equation}
\label{sos1z}
\hbox{\vtop{\offinterlineskip\halign{
\hfil#\hfil\cr
{\rm lim}\cr
$\stackrel{}{{}_{p_1,\ldots,p_k\to \infty}}$\cr
}} }\Biggl\Vert
K(t_1,\ldots,t_k)-
\sum_{j_1=0}^{p_1}\ldots \sum_{j_k=0}^{p_k}
C_{j_k\ldots j_1}\prod_{l=1}^{k} \phi_{j_l}(t_l)\Biggr
\Vert_{L_2([t, T]^k)}=0,
\end{equation}

\vspace{4mm}
\noindent
where

\vspace{-2mm}
\begin{equation}
\label{ppppa}
C_{j_k\ldots j_1}=\int\limits_{[t,T]^k}
K(t_1,\ldots,t_k)\prod_{l=1}^{k}\phi_{j_l}(t_l)dt_1\ldots dt_k
\end{equation}

\vspace{4mm}
\noindent
is the Fourier coefficient and

$$
\left\Vert f\right\Vert_{L_2([t, T]^k)}=\left(\int\limits_{[t,T]^k}
f^2(t_1,\ldots,t_k)dt_1\ldots dt_k\right)^{1/2}.
$$

\vspace{5mm}

Consider the partition $\{\tau_j\}_{j=0}^N$ of $[t,T]$ such that

\begin{equation}
\label{1111}
t=\tau_0<\ldots <\tau_N=T,\ \ \ \
\Delta_N=
\hbox{\vtop{\offinterlineskip\halign{
\hfil#\hfil\cr
{\rm max}\cr
$\stackrel{}{{}_{0\le j\le N-1}}$\cr
}} }\Delta\tau_j\to 0\ \ \hbox{if}\ \ N\to \infty,\ \ \ \
\Delta\tau_j=\tau_{j+1}-\tau_j.
\end{equation}

\vspace{2mm}

{\bf Theorem 1}\ \cite{2006} (2006), 
\cite{2011-2}-\cite{Mikh-1}. {\it Suppose that
every $\psi_l(\tau)$ $(l=1,\ldots, k)$ is a continuous nonrandom function on 
$[t, T]$ and
$\{\phi_j(x)\}_{j=0}^{\infty}$ is a complete orthonormal system  
of continous functions in $L_2([t,T]).$ Then

\vspace{1mm}
$$
J[\psi^{(k)}]_{T,t}=
\hbox{\vtop{\offinterlineskip\halign{
\hfil#\hfil\cr
{\rm l.i.m.}\cr
$\stackrel{}{{}_{p_1,\ldots,p_k\to \infty}}$\cr
}} }\sum_{j_1=0}^{p_1}\ldots\sum_{j_k=0}^{p_k}
C_{j_k\ldots j_1}\Biggl(
\prod_{l=1}^k\zeta_{j_l}^{(i_l)}-
\Biggr.
$$

\vspace{2mm}
\begin{equation}
\label{tyyy}
-\Biggl.
\hbox{\vtop{\offinterlineskip\halign{
\hfil#\hfil\cr
{\rm l.i.m.}\cr
$\stackrel{}{{}_{N\to \infty}}$\cr
}} }\sum_{(l_1,\ldots,l_k)\in {\rm G}_k}
\phi_{j_{1}}(\tau_{l_1})
\Delta{\bf w}_{\tau_{l_1}}^{(i_1)}\ldots
\phi_{j_{k}}(\tau_{l_k})
\Delta{\bf w}_{\tau_{l_k}}^{(i_k)}\Biggr),
\end{equation}

\vspace{5mm}
\noindent
where

\vspace{-2mm}
$$
{\rm G}_k={\rm H}_k\backslash{\rm L}_k,\ \ \ 
{\rm H}_k=\{(l_1,\ldots,l_k):\ l_1,\ldots,l_k=0,\ 1,\ldots,N-1\},
$$

$$
{\rm L}_k=\{(l_1,\ldots,l_k):\ l_1,\ldots,l_k=0,\ 1,\ldots,N-1;\
l_g\ne l_r\ (g\ne r);\ g, r=1,\ldots,k\},
$$

\vspace{4mm}
\noindent
${\rm l.i.m.}$ is a limit in the mean-square sense,\ 
$i_1,\ldots,i_k=0,1,\ldots,m,$\ 
$C_{j_k\ldots j_1}$ is the Fourier coefficient 
{\rm (\ref{ppppa}),}
\begin{equation}
\label{rr23}
\zeta_{j}^{(i)}=
\int\limits_t^T \phi_{j}(s) d{\bf w}_s^{(i)}
\end{equation} 

\vspace{2mm}
\noindent
are independent standard Gaussian random variables
for various
$i$ or $j$ {\rm(}in the case when $i\ne 0${\rm),}
$\Delta{\bf w}_{\tau_{j}}^{(i)}=
{\bf w}_{\tau_{j+1}}^{(i)}-{\bf w}_{\tau_{j}}^{(i)}$
$(i=0, 1,\ldots,m),$\
$\left\{\tau_{j}\right\}_{j=0}^{N}$ is a partition of
$[t,T],$ which satisfies the condition {\rm (\ref{1111})}.
}

\vspace{2mm}

In order to evaluate the significance of Theorem 1 for practice we will
demonstrate its transformed particular cases for 
$k=1,\ldots,5$ \cite{2006}-\cite{Mikh-1}

\vspace{1mm}
\begin{equation}
\label{a1}
J[\psi^{(1)}]_{T,t}
=\hbox{\vtop{\offinterlineskip\halign{
\hfil#\hfil\cr
{\rm l.i.m.}\cr
$\stackrel{}{{}_{p_1\to \infty}}$\cr
}} }\sum_{j_1=0}^{p_1}
C_{j_1}\zeta_{j_1}^{(i_1)},
\end{equation}

\vspace{2mm}

\begin{equation}
\label{a2}
J[\psi^{(2)}]_{T,t}
=\hbox{\vtop{\offinterlineskip\halign{
\hfil#\hfil\cr
{\rm l.i.m.}\cr
$\stackrel{}{{}_{p_1,p_2\to \infty}}$\cr
}} }\sum_{j_1=0}^{p_1}\sum_{j_2=0}^{p_2}
C_{j_2j_1}\Biggl(\zeta_{j_1}^{(i_1)}\zeta_{j_2}^{(i_2)}
-{\bf 1}_{\{i_1=i_2\ne 0\}}
{\bf 1}_{\{j_1=j_2\}}\Biggr),
\end{equation}

\vspace{4mm}

$$
J[\psi^{(3)}]_{T,t}=
\hbox{\vtop{\offinterlineskip\halign{
\hfil#\hfil\cr
{\rm l.i.m.}\cr
$\stackrel{}{{}_{p_1,\ldots,p_3\to \infty}}$\cr
}} }\sum_{j_1=0}^{p_1}\sum_{j_2=0}^{p_2}\sum_{j_3=0}^{p_3}
C_{j_3j_2j_1}\Biggl(
\zeta_{j_1}^{(i_1)}\zeta_{j_2}^{(i_2)}\zeta_{j_3}^{(i_3)}
-\Biggr.
$$

\begin{equation}
\label{a3}
-\Biggl.
{\bf 1}_{\{i_1=i_2\ne 0\}}
{\bf 1}_{\{j_1=j_2\}}
\zeta_{j_3}^{(i_3)}
-{\bf 1}_{\{i_2=i_3\ne 0\}}
{\bf 1}_{\{j_2=j_3\}}
\zeta_{j_1}^{(i_1)}-
{\bf 1}_{\{i_1=i_3\ne 0\}}
{\bf 1}_{\{j_1=j_3\}}
\zeta_{j_2}^{(i_2)}\Biggr),
\end{equation}

\vspace{4mm}

$$
J[\psi^{(4)}]_{T,t}
=
\hbox{\vtop{\offinterlineskip\halign{
\hfil#\hfil\cr
{\rm l.i.m.}\cr
$\stackrel{}{{}_{p_1,\ldots,p_4\to \infty}}$\cr
}} }\sum_{j_1=0}^{p_1}\ldots\sum_{j_4=0}^{p_4}
C_{j_4\ldots j_1}\Biggl(
\prod_{l=1}^4\zeta_{j_l}^{(i_l)}
\Biggr.
-
$$
$$
-
{\bf 1}_{\{i_1=i_2\ne 0\}}
{\bf 1}_{\{j_1=j_2\}}
\zeta_{j_3}^{(i_3)}
\zeta_{j_4}^{(i_4)}
-
{\bf 1}_{\{i_1=i_3\ne 0\}}
{\bf 1}_{\{j_1=j_3\}}
\zeta_{j_2}^{(i_2)}
\zeta_{j_4}^{(i_4)}-
$$
$$
-
{\bf 1}_{\{i_1=i_4\ne 0\}}
{\bf 1}_{\{j_1=j_4\}}
\zeta_{j_2}^{(i_2)}
\zeta_{j_3}^{(i_3)}
-
{\bf 1}_{\{i_2=i_3\ne 0\}}
{\bf 1}_{\{j_2=j_3\}}
\zeta_{j_1}^{(i_1)}
\zeta_{j_4}^{(i_4)}-
$$
$$
-
{\bf 1}_{\{i_2=i_4\ne 0\}}
{\bf 1}_{\{j_2=j_4\}}
\zeta_{j_1}^{(i_1)}
\zeta_{j_3}^{(i_3)}
-
{\bf 1}_{\{i_3=i_4\ne 0\}}
{\bf 1}_{\{j_3=j_4\}}
\zeta_{j_1}^{(i_1)}
\zeta_{j_2}^{(i_2)}+
$$
$$
+
{\bf 1}_{\{i_1=i_2\ne 0\}}
{\bf 1}_{\{j_1=j_2\}}
{\bf 1}_{\{i_3=i_4\ne 0\}}
{\bf 1}_{\{j_3=j_4\}}
+
$$
$$
+
{\bf 1}_{\{i_1=i_3\ne 0\}}
{\bf 1}_{\{j_1=j_3\}}
{\bf 1}_{\{i_2=i_4\ne 0\}}
{\bf 1}_{\{j_2=j_4\}}+
$$
\begin{equation}
\label{a4}
+\Biggl.
{\bf 1}_{\{i_1=i_4\ne 0\}}
{\bf 1}_{\{j_1=j_4\}}
{\bf 1}_{\{i_2=i_3\ne 0\}}
{\bf 1}_{\{j_2=j_3\}}\Biggr),
\end{equation}

\vspace{5mm}

$$
J[\psi^{(5)}]_{T,t}
=\hbox{\vtop{\offinterlineskip\halign{
\hfil#\hfil\cr
{\rm l.i.m.}\cr
$\stackrel{}{{}_{p_1,\ldots,p_5\to \infty}}$\cr
}} }\sum_{j_1=0}^{p_1}\ldots\sum_{j_5=0}^{p_5}
C_{j_5\ldots j_1}\Biggl(
\prod_{l=1}^5\zeta_{j_l}^{(i_l)}
-\Biggr.
$$
$$
-
{\bf 1}_{\{i_1=i_2\ne 0\}}
{\bf 1}_{\{j_1=j_2\}}
\zeta_{j_3}^{(i_3)}
\zeta_{j_4}^{(i_4)}
\zeta_{j_5}^{(i_5)}-
{\bf 1}_{\{i_1=i_3\ne 0\}}
{\bf 1}_{\{j_1=j_3\}}
\zeta_{j_2}^{(i_2)}
\zeta_{j_4}^{(i_4)}
\zeta_{j_5}^{(i_5)}-
$$
$$
-
{\bf 1}_{\{i_1=i_4\ne 0\}}
{\bf 1}_{\{j_1=j_4\}}
\zeta_{j_2}^{(i_2)}
\zeta_{j_3}^{(i_3)}
\zeta_{j_5}^{(i_5)}-
{\bf 1}_{\{i_1=i_5\ne 0\}}
{\bf 1}_{\{j_1=j_5\}}
\zeta_{j_2}^{(i_2)}
\zeta_{j_3}^{(i_3)}
\zeta_{j_4}^{(i_4)}-
$$
$$
-
{\bf 1}_{\{i_2=i_3\ne 0\}}
{\bf 1}_{\{j_2=j_3\}}
\zeta_{j_1}^{(i_1)}
\zeta_{j_4}^{(i_4)}
\zeta_{j_5}^{(i_5)}-
{\bf 1}_{\{i_2=i_4\ne 0\}}
{\bf 1}_{\{j_2=j_4\}}
\zeta_{j_1}^{(i_1)}
\zeta_{j_3}^{(i_3)}
\zeta_{j_5}^{(i_5)}-
$$
$$
-
{\bf 1}_{\{i_2=i_5\ne 0\}}
{\bf 1}_{\{j_2=j_5\}}
\zeta_{j_1}^{(i_1)}
\zeta_{j_3}^{(i_3)}
\zeta_{j_4}^{(i_4)}
-{\bf 1}_{\{i_3=i_4\ne 0\}}
{\bf 1}_{\{j_3=j_4\}}
\zeta_{j_1}^{(i_1)}
\zeta_{j_2}^{(i_2)}
\zeta_{j_5}^{(i_5)}-
$$
$$
-
{\bf 1}_{\{i_3=i_5\ne 0\}}
{\bf 1}_{\{j_3=j_5\}}
\zeta_{j_1}^{(i_1)}
\zeta_{j_2}^{(i_2)}
\zeta_{j_4}^{(i_4)}
-{\bf 1}_{\{i_4=i_5\ne 0\}}
{\bf 1}_{\{j_4=j_5\}}
\zeta_{j_1}^{(i_1)}
\zeta_{j_2}^{(i_2)}
\zeta_{j_3}^{(i_3)}+
$$
$$
+
{\bf 1}_{\{i_1=i_2\ne 0\}}
{\bf 1}_{\{j_1=j_2\}}
{\bf 1}_{\{i_3=i_4\ne 0\}}
{\bf 1}_{\{j_3=j_4\}}\zeta_{j_5}^{(i_5)}+
{\bf 1}_{\{i_1=i_2\ne 0\}}
{\bf 1}_{\{j_1=j_2\}}
{\bf 1}_{\{i_3=i_5\ne 0\}}
{\bf 1}_{\{j_3=j_5\}}\zeta_{j_4}^{(i_4)}+
$$
$$
+
{\bf 1}_{\{i_1=i_2\ne 0\}}
{\bf 1}_{\{j_1=j_2\}}
{\bf 1}_{\{i_4=i_5\ne 0\}}
{\bf 1}_{\{j_4=j_5\}}\zeta_{j_3}^{(i_3)}+
{\bf 1}_{\{i_1=i_3\ne 0\}}
{\bf 1}_{\{j_1=j_3\}}
{\bf 1}_{\{i_2=i_4\ne 0\}}
{\bf 1}_{\{j_2=j_4\}}\zeta_{j_5}^{(i_5)}+
$$
$$
+
{\bf 1}_{\{i_1=i_3\ne 0\}}
{\bf 1}_{\{j_1=j_3\}}
{\bf 1}_{\{i_2=i_5\ne 0\}}
{\bf 1}_{\{j_2=j_5\}}\zeta_{j_4}^{(i_4)}+
{\bf 1}_{\{i_1=i_3\ne 0\}}
{\bf 1}_{\{j_1=j_3\}}
{\bf 1}_{\{i_4=i_5\ne 0\}}
{\bf 1}_{\{j_4=j_5\}}\zeta_{j_2}^{(i_2)}+
$$
$$
+
{\bf 1}_{\{i_1=i_4\ne 0\}}
{\bf 1}_{\{j_1=j_4\}}
{\bf 1}_{\{i_2=i_3\ne 0\}}
{\bf 1}_{\{j_2=j_3\}}\zeta_{j_5}^{(i_5)}+
{\bf 1}_{\{i_1=i_4\ne 0\}}
{\bf 1}_{\{j_1=j_4\}}
{\bf 1}_{\{i_2=i_5\ne 0\}}
{\bf 1}_{\{j_2=j_5\}}\zeta_{j_3}^{(i_3)}+
$$
$$
+
{\bf 1}_{\{i_1=i_4\ne 0\}}
{\bf 1}_{\{j_1=j_4\}}
{\bf 1}_{\{i_3=i_5\ne 0\}}
{\bf 1}_{\{j_3=j_5\}}\zeta_{j_2}^{(i_2)}+
{\bf 1}_{\{i_1=i_5\ne 0\}}
{\bf 1}_{\{j_1=j_5\}}
{\bf 1}_{\{i_2=i_3\ne 0\}}
{\bf 1}_{\{j_2=j_3\}}\zeta_{j_4}^{(i_4)}+
$$
$$
+
{\bf 1}_{\{i_1=i_5\ne 0\}}
{\bf 1}_{\{j_1=j_5\}}
{\bf 1}_{\{i_2=i_4\ne 0\}}
{\bf 1}_{\{j_2=j_4\}}\zeta_{j_3}^{(i_3)}+
{\bf 1}_{\{i_1=i_5\ne 0\}}
{\bf 1}_{\{j_1=j_5\}}
{\bf 1}_{\{i_3=i_4\ne 0\}}
{\bf 1}_{\{j_3=j_4\}}\zeta_{j_2}^{(i_2)}+
$$
$$
+
{\bf 1}_{\{i_2=i_3\ne 0\}}
{\bf 1}_{\{j_2=j_3\}}
{\bf 1}_{\{i_4=i_5\ne 0\}}
{\bf 1}_{\{j_4=j_5\}}\zeta_{j_1}^{(i_1)}+
{\bf 1}_{\{i_2=i_4\ne 0\}}
{\bf 1}_{\{j_2=j_4\}}
{\bf 1}_{\{i_3=i_5\ne 0\}}
{\bf 1}_{\{j_3=j_5\}}\zeta_{j_1}^{(i_1)}+
$$
\begin{equation}
\label{a5}
+\Biggl.
{\bf 1}_{\{i_2=i_5\ne 0\}}
{\bf 1}_{\{j_2=j_5\}}
{\bf 1}_{\{i_3=i_4\ne 0\}}
{\bf 1}_{\{j_3=j_4\}}\zeta_{j_1}^{(i_1)}\Biggr),
\end{equation}

\vspace{5mm}
\noindent
where ${\bf 1}_A$ is the indicator of the set $A$.

It was shown 
that Theorem 1 is valid for convergence 
in the mean of degree $2n$ ($n\in \mathbb{N}$) 
\cite{2018a} (Sect.~1.1.9, 1.11, 1.12), \cite{arxiv-1} (Sect.~6, 15, 16)
and for convergence with probablity 1 (w.~p.~1) \cite{2018a} (Sect.~1.7.2), \cite{29}, \cite{new-new-2}
(the cases of Legendre polynomials and trigonometric 
functions).
Moreover, the complete orthonormal systems of Haar and 
Rademacher--Walsh functions in the space $L_2([t,T])$
can also be applied in Theorem 1
\cite{2006}-\cite{2013}.
The modification of Theorem 1 for 
complete orthonormal with weigth $r(x)\ge 0$ systems
of functions in the space $L_2([t,T])$ can be found in 
\cite{2018}-\cite{2018aaaaa1}, \cite{arxiv-26b}.
Recently, Theorem 1 and Theorem 2 (see below) has been applied 
to the expansion and mean-square
approximation of iterated stochastic integrals 
with respect to the infinite-dmensional $Q$-Wiener process
\cite{2018a}-\cite{2018aaaaa1} (Chapter 7), \cite{31a}, \cite{200a},  \cite{arxiv-21}-\cite{OK}. 
These results can be directly
applied to construction of 
high-order strong numerical methods for non-commutative 
semilinear stochastic partial differential equations
with multiplicative trace class noise
\cite{2018a}-\cite{2018aaaaa1} (Chapter 7), \cite{200a}, \cite{OK}.

Note that we obtain the following useful possibilities
of the approach based on Theorem 1.

1. There is the explicit formula (\ref{ppppa}) for calculation 
of expansion coefficients 
of the iterated Ito stochastic integral (\ref{ito}) with any
fixed multiplicity $k$. 

2. We have new possibilities for exact calculation of the mean-square 
approximation error
of iterated Ito stochastic integral (\ref{ito})
\cite{2017}-\cite{2018aaaaa1}, \cite{17a}, \cite{arxiv-2} (also see Theorem 3 below).

3. Since the used
multiple Fourier series is a generalized in the sense          
that it is built using various complete orthonormal
systems of functions in the space $L_2([t, T])$, then we 
have new possibilities 
for approximation --- we can
use not only trigonometric functions as in \cite{2}-\cite{3}
but Legendre polynomials.

4. As it turned out \cite{2006}-\cite{arxiv-9} it is more convenient to work 
with Legendre po\-ly\-no\-mi\-als for ap\-pro\-xi\-ma\-ti\-on
of the iterated Ito stochastic integrals (\ref{ito}). 
Approximations based on the Legendre polynomials are much simpler 
than their analogues based on the tri\-go\-no\-met\-ric functions.
Another advantages of the application of Legendre polynomials 
in the framework of the mentioned problem are considered
in \cite{2018a}-\cite{2018aaaaa1}, \cite{29a}, \cite{301a}.

5. The approach to expansion of iterated 
Ito and Stratonovich stochastic integrals 
based on the Karhunen--Loeve expansion
of the Brownian bridge process 
\cite{Mi2} (also see \cite{2}, \cite{3})
as well as the approach from \cite{rr}
lead to 
iterated application of the operation of limit
transition (the operation of limit transition 
is implemented only once in Theorem 1 and Theorem 2 (see below))
starting from  
the second multiplicity (in the general case) 
and third multiplicity (for the case
$\psi_1(s), \psi_2(s), \psi_3(s)\equiv 1;$ 
$i_1, i_2, i_3=1,\ldots,m$)
of iterated stochastic integrals.
Multiple series from Theorems 1, 2 (the operation of limit transition 
is implemented only once) are more convenient 
for approximation than the iterated ones
(iterated application of the operation of limit
transition), 
since partial sums of multiple series converge for any possible case of  
convergence to infinity of their upper limits of summation 
(let us denote them as $p_1,\ldots, p_k$). 
For example, when
$p_1=\ldots=p_k=p\to\infty$. 
For iterated series, the condition $p_1=\ldots=p_k=p\to\infty$ obviously 
does not guarantee the convergence of this series.
However, the authors of the works
\cite{2}
(Sect.~5.8, pp.~202--204), \cite{KPS} (pp.~82-84),
\cite{Zapad-2} (pp.~438-439),  
\cite{Zapad-9} (pp.~263-264) use 
the Wong--Zakai approximation 
\cite{W-Z-1}-\cite{Watanabe} (without rigorous proof) within the frames
of the method of expansion of iterated stochastic integrals
\cite{Mi2} (1988) based on the series expansion 
of the Brownian bridge process (version
of the so-called Karhunen-Loeve expansion).
See discussions in \cite{2018a} (Sect.~2.16, 6.2), 
\cite{2018aaaaa} (Sect.~2.6.2, 6.2),
\cite{2018aaaaa1} (Sect.~2.6.2, 6.2),
\cite{arxiv-1} (Sect.~11), \cite{arxiv-3} (Sect.~8),
\cite{arxiv-5} (Sect.~6) for detail.

Note that the correctness of the formulas (\ref{a1})--(\ref{a5}) 
can be 
verified 
by the fact that if 
$i_1=\ldots=i_5=i=1,\ldots,m$
and $\psi_1(s),\ldots,\psi_5(s)\equiv \psi(s)$ in (\ref{a1})--(\ref{a5}),
then we can obtain from (\ref{a1})--(\ref{a5}) the 
following equalities 

\vspace{-1mm}
$$
J[\psi^{(1)}]_{T,t}
=\frac{1}{1!}\delta_{T,t},
$$

$$
J[\psi^{(2)}]_{T,t}
=\frac{1}{2!}\left(\delta^2_{T,t}-\Delta_{T,t}\right),\
$$

$$
J[\psi^{(3)}]_{T,t}
=\frac{1}{3!}\left(\delta_{T,t}^3-3\delta_{T,t}\Delta_{T,t}\right),
$$

$$
J[\psi^{(4)}]_{T,t}
=\frac{1}{4!}\left(\delta^4_{T,t}-6\delta_{T,t}^2\Delta_{T,t}
+3\Delta^2_{T,t}\right),\
$$

$$
J[\psi^{(5)}]_{T,t}
=\frac{1}{5!}\left(\delta^5_{T,t}-10\delta_{T,t}^3\Delta_{T,t}
+15\delta_{T,t}\Delta^2_{T,t}\right)
$$

\vspace{4mm}
\noindent
w.~p.~1, where 

\vspace{-2mm}
$$
\delta_{T,t}=\int\limits_t^T\psi(s)d{\bf w}_s^{(i)},\ \ \
\Delta_{T,t}=\int\limits_t^T\psi^2(s)ds.
$$

\vspace{3mm}

The above formulas can be independently  
obtained using the Ito formula and Hermite polynomials.
Note that the cases  
$k=2, 3$ and $p_1=p_2=p_3=p$ are considered in detail
in \cite{2011-2}-\cite{2013}, \cite{arxiv-1}.

Consider the generalization of formulas (\ref{a1})--(\ref{a5}) 
for the case 
of an arbitrary multiplicity $k$ of the stochastic integral
$J[\psi^{(k)}]_{T,t}$ as well as 
for the case
of an arbitrary complete orthonormal systems  
of functions in the space $L_2([t,T])$ 
and $\psi_1(\tau),\ldots,\psi_k(\tau)\in L_2([t, T]).$
In order to do this, let us
consider the unordered
set $\{1, 2, \ldots, k\}$ 
and separate it into two parts:
the first part consists of $r$ unordered 
pairs (sequence order of these pairs is also unimportant) and the 
second one consists of the 
remaining $k-2r$ numbers.
So, we have

\vspace{-2mm}
\begin{equation}
\label{leto5007}
(\{
\underbrace{\{g_1, g_2\}, \ldots, 
\{g_{2r-1}, g_{2r}\}}_{\small{\hbox{part 1}}}
\},
\{\underbrace{q_1, \ldots, q_{k-2r}}_{\small{\hbox{part 2}}}
\}),
\end{equation}

\vspace{3mm}
\noindent
where 
$\{g_1, g_2, \ldots, 
g_{2r-1}, g_{2r}, q_1, \ldots, q_{k-2r}\}=\{1, 2, \ldots, k\},$
braces   
mean an unordered 
set, and pa\-ren\-the\-ses mean an ordered set.

We will say that (\ref{leto5007}) is a partition 
and consider the sum with respect to all possible
partitions

\begin{equation}
\label{leto5008}
\sum_{\stackrel{(\{\{g_1, g_2\}, \ldots, 
\{g_{2r-1}, g_{2r}\}\}, \{q_1, \ldots, q_{k-2r}\})}
{{}_{\{g_1, g_2, \ldots, 
g_{2r-1}, g_{2r}, q_1, \ldots, q_{k-2r}\}=\{1, 2, \ldots, k\}}}}
a_{g_1 g_2, \ldots, 
g_{2r-1} g_{2r}, q_1 \ldots q_{k-2r}}.
\end{equation}

\vspace{4mm}

Below there are several examples of sums in the form (\ref{leto5008})

\vspace{1mm}
$$
\sum_{\stackrel{(\{g_1, g_2\})}{{}_{\{g_1, g_2\}=\{1, 2\}}}}
a_{g_1 g_2}=a_{12},
$$

\vspace{3mm}
$$
\sum_{\stackrel{(\{\{g_1, g_2\}, \{g_3, g_4\}\})}
{{}_{\{g_1, g_2, g_3, g_4\}=\{1, 2, 3, 4\}}}}
a_{g_1 g_2 g_3 g_4}=a_{1234} + a_{1324} + a_{2314},
$$

\vspace{3mm}
$$
\sum_{\stackrel{(\{g_1, g_2\}, \{q_1, q_{2}\})}
{{}_{\{g_1, g_2, q_1, q_{2}\}=\{1, 2, 3, 4\}}}}
a_{g_1 g_2, q_1 q_{2}}=a_{12,34}+a_{13,24}+a_{14,23}
+a_{23,14}+a_{24,13}+a_{34,12},
$$

\vspace{3mm}
$$
\sum_{\stackrel{(\{g_1, g_2\}, \{q_1, q_{2}, q_3\})}
{{}_{\{g_1, g_2, q_1, q_{2}, q_3\}=\{1, 2, 3, 4, 5\}}}}
a_{g_1 g_2, q_1 q_{2}q_3}
=a_{12,345}+a_{13,245}+a_{14,235}
+a_{15,234}+a_{23,145}+a_{24,135}+
$$
$$
+a_{25,134}+a_{34,125}+a_{35,124}+a_{45,123},
$$

\vspace{3mm}
$$
\sum_{\stackrel{(\{\{g_1, g_2\}, \{g_3, g_{4}\}\}, \{q_1\})}
{{}_{\{g_1, g_2, g_3, g_{4}, q_1\}=\{1, 2, 3, 4, 5\}}}}
a_{g_1 g_2, g_3 g_{4},q_1}
=
a_{12,34,5}+a_{13,24,5}+a_{14,23,5}+
a_{12,35,4}+a_{13,25,4}+a_{15,23,4}+
$$
$$
+a_{12,54,3}+a_{15,24,3}+a_{14,25,3}+a_{15,34,2}+a_{13,54,2}+a_{14,53,2}+
a_{52,34,1}+a_{53,24,1}+a_{54,23,1}.
$$

\vspace{6mm}

Now we can generalize Theorem 1.

\vspace{2mm}

{\bf Theorem 2}\ \cite{2018a} (Sect.~1.11), \cite{arxiv-1} (Sect.~15).
{\it Suppose that
$\psi_1(\tau),\ldots,\psi_k(\tau)\in L_2([t, T])$ and
$\{\phi_j(x)\}_{j=0}^{\infty}$ is an arbitrary complete orthonormal system  
of functions in the space $L_2([t,T]).$
Then the following expansion

\vspace{1mm}
$$
J[\psi^{(k)}]_{T,t}=
\hbox{\vtop{\offinterlineskip\halign{
\hfil#\hfil\cr
{\rm l.i.m.}\cr
$\stackrel{}{{}_{p_1,\ldots,p_k\to \infty}}$\cr
}} }
\sum\limits_{j_1=0}^{p_1}\ldots
\sum\limits_{j_k=0}^{p_k}
C_{j_k\ldots j_1}\Biggl(
\prod_{l=1}^k\zeta_{j_l}^{(i_l)}+\sum\limits_{r=1}^{[k/2]}
(-1)^r \times \Biggr.
$$

\vspace{2mm}
\begin{equation}
\label{leto6000}
\times
\sum_{\stackrel{(\{\{g_1, g_2\}, \ldots, 
\{g_{2r-1}, g_{2r}\}\}, \{q_1, \ldots, q_{k-2r}\})}
{{}_{\{g_1, g_2, \ldots, 
g_{2r-1}, g_{2r}, q_1, \ldots, q_{k-2r}\}=\{1,2, \ldots, k\}}}}
\prod\limits_{s=1}^r
{\bf 1}_{\{i_{g_{{}_{2s-1}}}=~i_{g_{{}_{2s}}}\ne 0\}}
\Biggl.{\bf 1}_{\{j_{g_{{}_{2s-1}}}=~j_{g_{{}_{2s}}}\}}
\prod_{l=1}^{k-2r}\zeta_{j_{q_l}}^{(i_{q_l})}\Biggr)
\end{equation}

\vspace{5mm}
\noindent 
that converges in the mean-square sense is valid$,$ where 
$[x]$ is an integer part of a real number $x;$
another notations are the same as in Theorem~{\rm 1.}
}

\vspace{2mm}

In particular from (\ref{leto6000}) for $k=5$ we obtain

\vspace{1mm}

$$
J[\psi^{(5)}]_{T,t}=
\hbox{\vtop{\offinterlineskip\halign{
\hfil#\hfil\cr
{\rm l.i.m.}\cr
$\stackrel{}{{}_{p_1,\ldots,p_5\to \infty}}$\cr
}} }
\sum\limits_{j_1=0}^{p_1}\ldots
\sum\limits_{j_5=0}^{p_5}
C_{j_5 \ldots j_1}
\Biggl(
\zeta_{j_1}^{(i_1)}\ldots \zeta_{j_5}^{(i_5)}
-\Biggr.
$$

\vspace{2mm}
$$
-\sum\limits_{\stackrel{(\{g_1, g_2\}, \{q_1, q_{2}, q_3\})}
{{}_{\{g_1, g_2, q_{1}, q_{2}, q_3\}=\{1, 2, 3, 4, 5\}}}}
{\bf 1}_{\{i_{g_{{}_{1}}}=~i_{g_{{}_{2}}}\ne 0\}}
{\bf 1}_{\{j_{g_{{}_{1}}}=~j_{g_{{}_{2}}}\}}
\prod_{l=1}^{3}\zeta_{j_{q_l}}^{(i_{q_l})}+
$$

\vspace{2mm}
$$
\Biggl.+
\sum_{\stackrel{(\{\{g_1, g_2\}, 
\{g_{3}, g_{4}\}\}, \{q_1\})}
{{}_{\{g_1, g_2, g_{3}, g_{4}, q_1\}=\{1, 2, 3, 4, 5\}}}}
{\bf 1}_{\{i_{g_{{}_{1}}}=~i_{g_{{}_{2}}}\ne 0\}}
{\bf 1}_{\{j_{g_{{}_{1}}}=~j_{g_{{}_{2}}}\}}
\Biggl.{\bf 1}_{\{i_{g_{{}_{3}}}=~i_{g_{{}_{4}}}\ne 0\}}
{\bf 1}_{\{j_{g_{{}_{3}}}=~j_{g_{{}_{4}}}\}}
\zeta_{j_{q_1}}^{(i_{q_1})}\Biggr).
$$

\vspace{5mm}

The last equality obviously agrees with
(\ref{a5}).

It should be noted that an analogue of Theorem 2 for multiple Ito 
stochastic integrals was considered 
in \cite{Rybakov1000}. 
Note that we use another notations in comparison with \cite{Rybakov1000}.
Moreover, the proof of an analogue of Theorem 2
from \cite{Rybakov1000} is somewhat different from the proof given in 
\cite{2018a} (Sect.~1.11), \cite{arxiv-1} (Sect.~15).

Let us denote

$$
E_k^{p_1,\ldots,p_k}\stackrel{{\rm def}}
{=}{\sf M}\left\{\left(J[\psi^{(k)}]_{T,t}-
J[\psi^{(k)}]_{T,t}^{p_1,\ldots,p_k}\right)^2\right\},\ \ \ 
E_k^p\stackrel{{\rm def}}{=}E_k^{p_1,\ldots,p_k}\biggl|_{p_1=\ldots=p_k=p}\biggr.,
$$

\vspace{2mm}
$$
I_k\stackrel{{\rm def}}{=}\left\Vert K\right\Vert^2_{L_2([t,T]^k)}=\int\limits_{[t,T]^k}
K^2(t_1,\ldots,t_k)dt_1\ldots dt_k,
$$

\vspace{4mm}
\noindent
where $J[\psi^{(k)}]_{T,t}^{p_1,\ldots,p_k}$ is the expression
on the right-hand side of (\ref{leto6000}) before passing to the limit
$\hbox{\vtop{\offinterlineskip\halign{
\hfil#\hfil\cr
{\rm l.i.m.}\cr
$\stackrel{}{{}_{p_1,\ldots,p_k\to \infty}}$\cr
}} }$, i.e.

$$
J[\psi^{(k)}]_{T,t}^{p_1,\ldots,p_k}=
\sum\limits_{j_1=0}^{p_1}\ldots
\sum\limits_{j_k=0}^{p_k}
C_{j_k\ldots j_1}\Biggl(
\prod_{l=1}^k\zeta_{j_l}^{(i_l)}+\sum\limits_{r=1}^{[k/2]}
(-1)^r \times \Biggr.
$$

\vspace{3mm}
$$
\times
\sum_{\stackrel{(\{\{g_1, g_2\}, \ldots, 
\{g_{2r-1}, g_{2r}\}\}, \{q_1, \ldots, q_{k-2r}\})}
{{}_{\{g_1, g_2, \ldots, 
g_{2r-1}, g_{2r}, q_1, \ldots, q_{k-2r}\}=\{1,2, \ldots, k\}}}}
\prod\limits_{s=1}^r
{\bf 1}_{\{i_{g_{{}_{2s-1}}}=~i_{g_{{}_{2s}}}\ne 0\}}
\Biggl.{\bf 1}_{\{j_{g_{{}_{2s-1}}}=~j_{g_{{}_{2s}}}\}}
\prod_{l=1}^{k-2r}\zeta_{j_{q_l}}^{(i_{q_l})}\Biggr).
$$

\vspace{5mm}

In \cite{2018a}-\cite{2018aaaaa1}, \cite{arxiv-1} it was shown that 

$$
E_k^{p_1,\ldots,p_k}\le k!\left(I_k-\sum_{j_1=0}^{p_1}\ldots
\sum_{j_k=0}^{p_k}C^2_{j_k\ldots j_1}\right)
$$

\vspace{4mm}
\noindent
if $i_1,\ldots,i_k=1,\ldots,m$ and $0<T-t<\infty$ or 
$i_1,\ldots,i_k=0, 1,\ldots,m$ and $0<T-t<1.$

Moreover, 
in \cite{2018a} (Sect.~1.1.9, 1.11, 1.12), \cite{arxiv-1} (Sect.~6, 15, 16)
it was shown that

\vspace{2mm}
$$
E_k^{p_1,\ldots,p_k}\le (k!)^{n} (2n-1)^{nk}\times
$$

\begin{equation}
\label{99999}
\times \left(I_k-\sum_{j_1=0}^{p_1}\ldots
\sum_{j_k=0}^{p_k}C^2_{j_k\ldots j_1}\right)^n,
\end{equation}

\vspace{4mm}
\noindent
where $n\in \mathbb{N}$.

The value $E_k^{p}$
can be calculated exactly.

\vspace{2mm}

{\bf Theorem 3} \cite{2018a} (Sect.~1.12), \cite{arxiv-2} (Sect.~6).
{\it Suppose that $\{\phi_j(x)\}_{j=0}^{\infty}$ 
is an arbitrary complete orthonormal system  
of functions in the space $L_2([t,T])$ and
$\psi_1(\tau),\ldots,\psi_k(\tau)\in L_2([t, T]).$  
Then

\begin{equation}
\label{tttr11}
E_k^p=I_k- \sum_{j_1,\ldots, j_k=0}^{p}
C_{j_k\ldots j_1}
{\sf M}\left\{J[\psi^{(k)}]_{T,t}
\sum\limits_{(j_1,\ldots,j_k)}
\int\limits_t^T \phi_{j_k}(t_k)
\ldots
\int\limits_t^{t_{2}}\phi_{j_{1}}(t_{1})
d{\bf f}_{t_1}^{(i_1)}\ldots
d{\bf f}_{t_k}^{(i_k)}\right\},
\end{equation}

\vspace{5mm}
\noindent
where $i_1,\ldots,i_k = 1,\ldots,m;$
the expression 

\vspace{-1mm}
$$
\sum\limits_{(j_1,\ldots,j_k)}
$$ 

\vspace{3mm}
\noindent
means the sum with respect to all
possible permutations 
$(j_1,\ldots,j_k)$. At the same time if 
$j_r$ swapped with $j_q$ in the permutation $(j_1,\ldots,j_k),$
then $i_r$ swapped with $i_q$ in the permutation
$(i_1,\ldots,i_k);$
another notations are the same as in Theorems {\rm 1, 2.}
}

\vspace{2mm}

Note that 

$$
{\sf M}\left\{J[\psi^{(k)}]_{T,t}
\int\limits_t^T \phi_{j_k}(t_k)
\ldots
\int\limits_t^{t_{2}}\phi_{j_{1}}(t_{1})
d{\bf f}_{t_1}^{(i_1)}\ldots
d{\bf f}_{t_k}^{(i_k)}\right\}=C_{j_k\ldots j_1}.
$$

\vspace{4mm}

Then from Theorem 3 for pairwise different $i_1,\ldots,i_k$ 
and for $i_1=\ldots=i_k$
we obtain

\vspace{2mm}
$$
E_k^p= I_k- \sum_{j_1,\ldots,j_k=0}^{p}
C_{j_k\ldots j_1}^2,
$$

\vspace{3mm}
$$ 
E_k^p= I_k - \sum_{j_1,\ldots,j_k=0}^{p}
C_{j_k\ldots j_1}\Biggl(\sum\limits_{(j_1,\ldots,j_k)}
C_{j_k\ldots j_1}\Biggr).
$$

\vspace{5mm}

\section{New Representation of the Levy Stochastic Area
Based on the Legendre Polynomials}

\vspace{5mm}

Let us consider (\ref{a2}) for the case 
$i_1\ne i_2$,\ \ $\psi_1(s),$ $\psi_2(s)\equiv 1$. At that we suppose that
$\{\phi_j(x)\}_{j=0}^{\infty}$ is the complete orthonormal sys\-tem  
of Legendre polynomials in the space $L_2([t,T]).$ Then 

\vspace{-1mm}
\begin{equation}
\label{zorro}
I_{T,t}^{(i_1 i_2)}
=
\frac{T-t}{2}\left(\zeta_0^{(i_1)}\zeta_0^{(i_2)}+\sum_{i=1}^{\infty}
\frac{1}{\sqrt{4i^2-1}}\left(
\zeta_{i-1}^{(i_1)}\zeta_{i}^{(i_2)}-
\zeta_i^{(i_1)}\zeta_{i-1}^{(i_2)}\right)\right),
\end{equation}

\vspace{3mm}
\noindent
where
\begin{equation}
\label{jq2}
I_{T,t}^{(i_1i_2)}=\int\limits_t^T\int\limits_t^s
d{\bf w}_{\tau}^{(i_1)}d{\bf w}_s^{(i_2)}\ \ \ (i_1,i_2=1,\ldots,m),
\end{equation}

\vspace{2mm}
\noindent
$\zeta_j^{(i)}$ 
are independent standard Gaussian random variables
(for various
$i$ or $j$), which have the following form  
$$
\zeta_j^{(i)}=\int\limits_t^T
\phi_j(s)d{\bf w}_s^{(i)},
$$

\vspace{1mm}
\noindent
where 
\begin{equation}
\label{fi}
\phi_i(s)=\sqrt{\frac{2i+1}{T-t}}P_i\left(\left(s-t-\frac{T-t}{2}\right)
\frac{2}{T-t}\right),\ \ \ i=0, 1, 2,\ldots, 
\end{equation}

\vspace{3mm}
\noindent
and $P_i(x)$ $(i=0, 1, 2,\ldots)$ is the Legendre polynomial.

Note that the representation (\ref{zorro}) was first
obtained in the author's works \cite{300a} (1997), \cite{400a} (1998).

From (\ref{zorro}) we obtain

\vspace{1mm}
$$
\frac{T-t}{2}\sum_{i=1}^{\infty}
\frac{1}{\sqrt{4i^2-1}}\left(
\zeta_{i-1}^{(i_1)}\zeta_{i}^{(i_2)}-
\zeta_i^{(i_1)}\zeta_{i-1}^{(i_2)}\right)
=\frac{1}{2}\left(I_{T,t}^{(i_1 i_2)}-I_{T,t}^{(i_2 i_1)}\right).
$$

\vspace{4mm}

Then,
a new representation of the Levy stochastic area based 
on the Legendre polynomials
has the following form
\begin{equation}
\label{l100}
A_{T,t}^{(i_1 i_2)}=
\frac{T-t}{2}\sum_{i=1}^{\infty}
\frac{1}{\sqrt{4i^2-1}}\left(
\zeta_{i-1}^{(i_1)}\zeta_{i}^{(i_2)}-
\zeta_i^{(i_1)}\zeta_{i-1}^{(i_2)}\right).
\end{equation}

\vspace{5mm}

\section{The Classical Representation of the Levy Stochastic Area}

\vspace{5mm}

Let us consider (\ref{a2}) for the case 
$i_1\ne i_2$,\ \ $\psi_1(s),$ $\psi_2(s)\equiv 1$. At that we suppose that
$\{\phi_j(x)\}_{j=0}^{\infty}$ is the complete orthonormal system  
of trigonometric functions in $L_2([t,T]).$ Then

\vspace{1mm}
$$
I_{T,t}^{(i_1 i_2)}=
\frac{1}{2}(T-t)\Biggl(
\zeta_{0}^{(i_1)}\zeta_{0}^{(i_2)}
+\frac{1}{\pi}
\sum_{r=1}^{\infty}\frac{1}{r}\biggl(
\zeta_{2r}^{(i_1)}\zeta_{2r-1}^{(i_2)}-
\zeta_{2r-1}^{(i_1)}\zeta_{2r}^{(i_2)}+\biggr.\Biggr.
$$

\vspace{1mm}
\begin{equation}
\label{ajjja}
+\biggl.\Biggl.
\sqrt{2}\left(\zeta_{2r-1}^{(i_1)}\zeta_{0}^{(i_2)}-
\zeta_{0}^{(i_1)}\zeta_{2r-1}^{(i_2)}\right)\biggr)\Biggr),
\end{equation}

\vspace{4mm}
\noindent
where we use the same notations as in (\ref{zorro}), but $\phi_j(s)$
has the following form

\vspace{1mm}

\begin{equation}
\label{rre}
\phi_j(s)=\frac{1}{\sqrt{T-t}}
\begin{cases}
1,\ &\ {\rm if}\ j=0\\
~\\
\sqrt{2}{\rm sin}(2\pi r(s-t)/(T-t)),\ &\ {\rm if}\ j=2r-1\\
~\\
\sqrt{2}{\rm cos}(2\pi r(s-t)/(T-t)),\ &\ {\rm if}\ j=2r
\end{cases},\ \ \ r=1,\ 2,\ldots
\end{equation}

\vspace{6mm}

From (\ref{ajjja}) we obtain

\vspace{1mm}
$$
\frac{T-t}{2\pi}
\sum_{r=1}^{\infty}\frac{1}{r}\left(
\zeta_{2r}^{(i_1)}\zeta_{2r-1}^{(i_2)}-
\zeta_{2r-1}^{(i_1)}\zeta_{2r}^{(i_2)}+
\sqrt{2}\left(\zeta_{2r-1}^{(i_1)}\zeta_{0}^{(i_2)}-
\zeta_{0}^{(i_1)}\zeta_{2r-1}^{(i_2)}\right)\right)
=
$$

\vspace{2mm}
$$
=\frac{1}{2}\left(I_{T,t}^{(i_1 i_2)}-I_{T,t}^{(i_2 i_1)}\right).
$$

\vspace{5mm}

Then, 
the representation of the Levy stochastic area based 
on the trigonometric functions
has the following form

\vspace{1mm}
\begin{equation}
\label{l200}
{\hat A}_{T,t}^{(i_1 i_2)}=
\frac{T-t}{2\pi}
\sum_{r=1}^{\infty}\frac{1}{r}\left(
\zeta_{2r}^{(i_1)}\zeta_{2r-1}^{(i_2)}-
\zeta_{2r-1}^{(i_1)}\zeta_{2r}^{(i_2)}+
\sqrt{2}\left(\zeta_{2r-1}^{(i_1)}\zeta_{0}^{(i_2)}-
\zeta_{0}^{(i_1)}\zeta_{2r-1}^{(i_2)}\right)\right).
\end{equation}

\vspace{4mm}

As we mentioned above,
Milstein G.N. proposed \cite{Mi2} the method of expansion 
of iterated Ito stochastic integrals of multiplicity 2
based on the trigonometric Fourier expansion of the 
following Brownian bridge 
process 

\vspace{-1mm}
$$
{\bf w}_t-\frac{t}{\Delta}{\bf w}_{\Delta},\ \ \
t\in[0,\Delta],\ \ \ \Delta>0,
$$

\vspace{3mm}
\noindent
where ${\bf w}_t$ is a standard
multidimensional Wiener process with independent components
${\bf w}^{(i)}_t$,\ $i=1,\ldots,m.$

The trigonometric Fourier expansion of the Brownian bridge 
process (version of the so-called Karunen--Loeve expansion) has the form
\cite{Mi2}

\begin{equation}
\label{6.5.2}
{\bf w}_t^{(i)}-\frac{t}{\Delta}{\bf w}_{\Delta}^{(i)}=
\frac{1}{2}a_{i,0}+\sum_{r=1}^{\infty}\left(
a_{i,r}{\rm cos}\frac{2\pi rt}{\Delta} +b_{i,r}{\rm sin}
\frac{2\pi rt}{\Delta}\right),
\end{equation}

\vspace{3mm}
\noindent
where

\vspace{-1mm}
$$
a_{i,r}=\frac{2}{\Delta} \int\limits_0^{\Delta}
\left({\bf w}_s^{(i)}-\frac{s}{\Delta}{\bf w}_{\Delta}^{(i)}\right)
{\rm cos}\frac{2\pi rs}{\Delta}ds,\
$$
$$
b_{i,r}=\frac{2}{\Delta} \int\limits_0^{\Delta}
\left({\bf w}_s^{(i)}-\frac{s}{\Delta}{\bf w}_{\Delta}^{(i)}\right)
{\rm sin}\frac{2\pi rs}{\Delta}ds,
$$

\vspace{3mm}
\noindent
$r=0, 1,\ldots,$\ \ $i=1,\ldots,m.$ 

It is easy to demonstrate \cite{Mi2} that the random variables
$a_{i,r}, b_{i,r}$ 
are Gaussian ones and they satisfy the following relations

\vspace{-2mm}
$$
{\sf M}\left\{a_{i,r}b_{i,r}\right\}=
{\sf M}\left\{a_{i,r}b_{i,k}\right\}=0,
$$

$$
{\sf M}\left\{a_{i,r}a_{i,k}\right\}=
{\sf M}\left\{b_{i,r}b_{i,k}\right\}=0,
$$

$$
{\sf M}\left\{a_{i_1,r}a_{i_2,r}\right\}=
{\sf M}\left\{b_{i_1,r}b_{i_2,r}\right\}=0,
$$

$$
{\sf M}\left\{a_{i,r}^2\right\}=
{\sf M}\left\{b_{i,r}^2\right\}=\frac{\Delta}{2\pi^2 r^2},
$$

\vspace{6mm}
\noindent
where $i,\ i_1, i_2=1,\ldots,m,$\ \ $r\ne k,$\ \ $i_1\ne i_2.$

According to (\ref{6.5.2}), we have

\begin{equation}
\label{6.5.7}
{\bf w}_t^{(i)}={\bf w}_{\Delta}^{(i)}\frac{t}{\Delta}+
\frac{1}{2}a_{i,0}+
\sum_{r=1}^{\infty}\left(
a_{i,r}{\rm cos}\frac{2\pi rt}{\Delta}+b_{i,r}{\rm sin}
\frac{2\pi rt}{\Delta}\right),
\end{equation}

\vspace{5mm}
\noindent
where the series
converges in the mean-square sense.

The expansion (\ref{ajjja}) has been obtained in \cite{Mi2} using 
(\ref{6.5.7}).

\vspace{5mm}

\section{New Simple Method for Obtainment of Representation of the Levy 
Stochastic Area}

\vspace{5mm}

It is well known that the idea of representing of the Wiener 
process as a functional series with random coefficients using the
complete orthonormal system of trigonometric functions in $L_2([0, T])$ 
goes back to the works of Wiener \cite{7b} (1924) and Levy 
\cite{7c} (1951). 
The specified series was used in \cite{7b} and \cite{7c} 
for construction of the Brownian motion process (Wiener process). 
A little later, Ito and McKean in \cite{7d} (1965) used for 
this purpose the complete orthonormal system of
Haar functions in $L_2([0, T])$.

Let ${\bf w}_{\tau},$ $\tau\in[0, T]$ be an $m$-dimestional
standard Wiener process with independent components
${\bf w}_{\tau}^{(i)}$ $(i=1,\ldots,m).$ 
We have

$$
{\bf w}_s^{(i)}-{\bf w}_t^{(i)}=
\int\limits_t^s d{\bf w}_{\tau}^{(i)}=
\int\limits_t^T {\bf 1}_{\{\tau<s\}} d{\bf w}_{\tau}^{(i)},
$$

\vspace{3mm}
\noindent
where
$$
{\bf 1}_{\{\tau<s\}}=
\begin{cases}
1,\ &\ \tau<s\\
~\\
0,\ &\ \hbox{\rm otherwise}
\end{cases},\ \ \ \tau, s\in [t, T],\ \ \ 0\le t<T.
$$

\vspace{5mm}

Consider the Fourier expansion of ${\bf 1}_{\{\tau<s\}}$ at the 
interval $[t, T]$ (see, for example, \cite{7e})

\begin{equation}
\label{eee}
{\bf 1}_{\{\tau<s\}}=\sum_{j=0}^{\infty}\int\limits_t^T
{\bf 1}_{\{\tau<s\}}\phi_j(\tau)d\tau \cdot \phi_j(\tau)=
\sum_{j=0}^{\infty}\int\limits_t^s
\phi_j(\tau)d\tau \cdot \phi_j(\tau),
\end{equation}

\vspace{4mm}
\noindent
where $\{\phi_j(\tau)\}_{j=0}^{\infty}$ is a complete 
orthonormal system of functions in the space $L_2([t, T])$
and the series on the right-hand side of (\ref{eee}) converges
in the mean-square sence, i.e.

$$
\int\limits_t^T\left({\bf 1}_{\{\tau<s\}}-
\sum_{j=0}^{q}\int\limits_t^s
\phi_j(\tau)d\tau \cdot \phi_j(\tau)\right)^2 d\tau\to 0\ \ \
\hbox{if}\ \ \ q\to\infty.
$$

\vspace{5mm}

Let $\left({\bf w}_s^{(i)}-{\bf w}_t^{(i)}\right)^{(q)}$ be 
the mean-square approximation of the process
${\bf w}_s^{(i)}-{\bf w}_t^{(i)}$, 
which has the following form

\begin{equation}
\label{jq1}
\left({\bf w}_s^{(i)}-{\bf w}_t^{(i)}\right)^{(q)}=
\int\limits_t^T
\left(\sum_{j=0}^{q}\int\limits_t^s
\phi_j(\tau)d\tau \cdot \phi_j(\tau)\right)
d{\bf w}_{\tau}^{(i)}=\sum_{j=0}^{q}\int\limits_t^s
\phi_j(\tau)d\tau   \cdot \int\limits_t^T\phi_j(\tau)
d{\bf w}_{\tau}^{(i)}.
\end{equation}

\vspace{4mm}

Moreover,

\vspace{-1mm}
$$
{\sf M}\left\{\Biggl(
{\bf w}_s^{(i)}-{\bf w}_t^{(i)}-
\left({\bf w}_s^{(i)}-{\bf w}_t^{(i)}\right)^{(q)}\Biggr)^2\right\}=
$$

$$
={\sf M}\left\{\left(
\int\limits_t^T\left(
{\bf 1}_{\{\tau<s\}}-
\sum_{j=0}^{q}\int\limits_t^s
\phi_j(\tau)d\tau \cdot \phi_j(\tau)\right)
d{\bf w}_{\tau}^{(i)}\right)^2\right\}=
$$

\begin{equation}
\label{t1}
=\int\limits_t^T
\left({\bf 1}_{\{\tau<s\}}-
\sum_{j=0}^{q}\int\limits_t^s
\phi_j(\tau)d\tau \cdot \phi_j(\tau)\right)^2 d\tau\to 0\ \ \
\hbox{if}\ q\to\infty.
\end{equation}

\vspace{5mm}

In \cite{rr} it was proposed to use the expansion similar 
to (\ref{jq1}) for construction of expansion of the iterated Ito stochastic 
integral (\ref{jq2}) of multiplicity 2. At that,
to obtain the mentioned expansion of (\ref{jq2}), the truncated 
expansions (\ref{jq1}) of components of the Wiener 
process ${\bf w}_s$ have been
iteratively substituted in the single integrals \cite{rr}. 
This procedure leads to the calculation
of coefficients of the double Fourier series, 
which is a time-consuming task for not too complex problem
of expansion of the iterated Ito stochastic integral (\ref{jq2}).

In contrast to \cite{rr} we subsitute the truncated expansion 
(\ref{jq1}) only one time and only into the innermost integral
in (\ref{jq2}). 
This procedure leads to the simple calculation
of the coefficients 

\vspace{-1mm}
$$
\int\limits_t^s
\phi_j(\tau)d\tau\ \ \ (j=0, 1, 2,\ldots)
$$

\vspace{3mm}
\noindent
of the usual (not double) Fourier series.

Moreover, we use the Legendre polynomials for construction 
of the expansion of (\ref{jq2}). For the first time the 
Legendre polynomials have been applied in the 
framework of the mentioned  problem in the author's papers 
\cite{300a} (1997), \cite{400a} (1998), 
\cite{500a} (2000), 
\cite{600a} (2001)
(also see \cite{2006}-\cite{Mikh-1}, \cite{new-art-1}, \cite{last-1}).
At the same time
in the papers of other author's these polynomials 
have not been considered as the basis functions for construction
of expansions of iterated Ito and Stratonovich 
stochastic integrals.

\vspace{2mm}

{\bf Theorem 4}\ \cite{2018a}-\cite{2018aaaaa1}, \cite{new-art-1}, \cite{last-1}. {\it Let
$\phi_j(\tau)$ $(j=0, 1, \ldots )$ be an arbitrary complete
orthonormal system of functions in the space $L_2([t, T]).$
Let 

\vspace{-1mm}
\begin{equation}
\label{l1}
\int\limits_t^T
\left({\bf w}_s^{(i_1)}-{\bf w}_t^{(i_1)}\right)^{(q)}
d{\bf w}_s^{(i_2)}=
\sum_{j=0}^{q}\int\limits_t^T\phi_j(\tau)
d{\bf w}_{\tau}^{(i_1)}
\int\limits_t^T\int\limits_t^s
\phi_j(\tau)d\tau d{\bf w}_{s}^{(i_2)}
\end{equation}

\vspace{3mm}
\noindent
be the approximation of the iterated Ito stochastic integral

$$
\int\limits_t^T
\int\limits_t^s
d{\bf w}_{\tau}^{(i_1)}d{\bf w}_{s}^{(i_2)}\ \ \ (i_1\ne i_2),
$$

\vspace{3mm}
\noindent
where $i_1,i_2=1,\ldots,m$.
Then 

$$
\int\limits_t^T
\int\limits_t^s
d{\bf w}_{\tau}^{(i_1)}d{\bf w}_{s}^{(i_2)}
=\hbox{\vtop{\offinterlineskip\halign{
\hfil#\hfil\cr
{\rm l.i.m.}\cr
$\stackrel{}{{}_{q\to \infty}}$\cr
}} }
\int\limits_t^T
\left({\bf w}_s^{(i_1)}-{\bf w}_t^{(i_1)}\right)^{(q)}
d{\bf w}_s^{(i_2)}=
$$

\vspace{2mm}
$$
=
\hbox{\vtop{\offinterlineskip\halign{
\hfil#\hfil\cr
{\rm l.i.m.}\cr
$\stackrel{}{{}_{q\to \infty}}$\cr
}} }
\sum_{j=0}^{q}\int\limits_t^T\phi_j(\tau)
d{\bf w}_{\tau}^{(i_1)}
\int\limits_t^T\int\limits_t^s
\phi_j(\tau)d\tau d{\bf w}_{s}^{(i_2)},
$$

\vspace{5mm}
\noindent
where $i_1\ne i_2$ $(i_1, i_2=1,\ldots,m)$.}

\vspace{2mm}

{\bf Proof.} Using standard properties of the Ito stochastic
integral as well as (\ref{t1}) and the property of orthonormality
of the functions 
$\phi_j(\tau)$ $(j=0, 1, \ldots )$ at the interval $[t, T],$
we obtain 

\vspace{1mm}
$$
{\sf M}\left\{\left(
\int\limits_t^T
\int\limits_t^s
d{\bf w}_{\tau}^{(i_1)}d{\bf w}_{s}^{(i_2)}-
\int\limits_t^T
\left({\bf w}_s^{(i_1)}-{\bf w}_t^{(i_1)}\right)^{(q)}
d{\bf w}_s^{(i_2)}\right)^2\right\}=
$$

\vspace{2mm}
$$
=
\int\limits_t^T
{\sf M}\left\{\Biggl(
{\bf w}_s^{(i_1)}-{\bf w}_t^{(i_1)}-
\left({\bf w}_s^{(i_1)}-{\bf w}_t^{(i_1)}\right)^{(q)}\Biggr)^2\right\}ds=
$$

\vspace{2mm}
$$
=
\int\limits_t^T\int\limits_t^T
\left({\bf 1}_{\{\tau<s\}}-
\sum_{j=0}^{q}\int\limits_t^s
\phi_j(\tau)d\tau \cdot \phi_j(\tau)\right)^2d\tau ds=
$$

\vspace{2mm}
\begin{equation}
\label{l5}
=\int\limits_t^T\left((s-t)-
\sum_{j=0}^{q}\left(\int\limits_t^s
\phi_j(\tau)d\tau\right)^2\right) ds.
\end{equation}

\vspace{5mm}

Applying the continuity of the functions $u_q(s)$ (see below),
the nondecreasing property of
the functional sequence

\vspace{-2mm}
$$
u_q(s)=\sum_{j=0}^{q}\left(\int\limits_t^s
\phi_j(\tau)d\tau\right)^2,
$$

\vspace{4mm}
\noindent
and the continuity of the limit function
$u(s)= s-t$
according to Dini's 
Theorem,
we have the uniform convergence 
$u_q(s)$ to $u(s)$ at the interval $[t, T]$.

Then from this fact as well as from (\ref{l5}) we obtain

\vspace{2mm}

\begin{equation}
\label{rez1}
\int\limits_t^T
\int\limits_t^s
d{\bf w}_{\tau}^{(i_1)}d{\bf w}_{s}^{(i_2)}
=\hbox{\vtop{\offinterlineskip\halign{
\hfil#\hfil\cr
{\rm l.i.m.}\cr
$\stackrel{}{{}_{q\to \infty}}$\cr
}} }
\int\limits_t^T
\left({\bf w}_s^{(i_1)}-{\bf w}_t^{(i_1)}\right)^{(q)}
d{\bf w}_s^{(i_2)}.
\end{equation}

\vspace{5mm}

Theorem 4 is proved.

Let $\{\phi_j(\tau)\}_{j=0}^{\infty}$ be the complete 
orthonormal system of Legendre polynomials in the space $L_2([t, T]),$
which has the form (\ref{fi}).
Then

\begin{equation}
\label{l2}
\int\limits_t^s
\phi_j(\tau)d\tau=
\frac{T-t}{2}\left(\frac{\phi_{j+1}(s)}{\sqrt{(2j+1)(2j+3)}}-
\frac{\phi_{j-1}(s)}{\sqrt{4j^2-1}}\right)\ \ \ \hbox{for}\ \ \ j\ge 1.
\end{equation}

\vspace{4mm}

Denote (see Theorem 1)
$$
\zeta_j^{(i)}=\int\limits_t^T\phi_j(\tau)
d{\bf w}_{\tau}^{(i)}\ \ \ (i=1,\ldots,m).
$$

\vspace{3mm}

From (\ref{l1}) and (\ref{l2}) we get

$$
\int\limits_t^T
\left({\bf w}_s^{(i_1)}-{\bf w}_t^{(i_1)}\right)^{(q)}
d{\bf w}_s^{(i_2)}=
\frac{1}{\sqrt{T-t}}\zeta_0^{(i_1)}
\int\limits_t^T(s-t){\bf w}_s^{(i_2)}+
$$

\vspace{1mm}
$$
+
\frac{T-t}{2}\sum_{j=1}^q
\zeta_j^{(i_1)}
\left(\frac{1}{\sqrt{(2j+1)(2j+3)}}\zeta_{j+1}^{(i_2)}-
\frac{1}{\sqrt{4j^2-1}}\zeta_{j-1}^{(i_2)}\right)=
$$

\vspace{4mm}
$$
=
\frac{T-t}{2}\zeta_0^{(i_1)}\left(
\zeta_0^{(i_2)}+\frac{1}{\sqrt{3}}\zeta_1^{(i_2)}\right)
+
$$

\vspace{1mm}
$$
+
\frac{T-t}{2}\sum_{j=1}^q
\zeta_j^{(i_1)}
\left(\frac{1}{\sqrt{(2j+1)(2j+3)}}\zeta_{j+1}^{(i_2)}-
\frac{1}{\sqrt{4j^2-1}}\zeta_{j-1}^{(i_2)}\right)=
$$

\vspace{4mm}
$$
=\frac{T-t}{2}\left(\zeta_0^{(i_1)}\zeta_0^{(i_2)}+\sum_{j=1}^{q}
\frac{1}{\sqrt{4j^2-1}}\left(
\zeta_{j-1}^{(i_1)}\zeta_{j}^{(i_2)}-
\zeta_j^{(i_1)}\zeta_{j-1}^{(i_2)}\right)\right)+
$$

\vspace{1mm}
\begin{equation}
\label{l6}
+
\frac{T-t}{2}\zeta_q^{(i_1)}\zeta_{q+1}^{(i_2)}
\frac{1}{\sqrt{(2q+1)(2q+3)}}.
\end{equation}

\vspace{5mm}

Then from (\ref{rez1}) and (\ref{l6}) we obtain

\vspace{1mm}
$$
\int\limits_t^T
\int\limits_t^s
d{\bf w}_{\tau}^{(i_1)}d{\bf w}_{s}^{(i_2)}
=\hbox{\vtop{\offinterlineskip\halign{
\hfil#\hfil\cr
{\rm l.i.m.}\cr
$\stackrel{}{{}_{q\to \infty}}$\cr
}} }
\int\limits_t^T
\left({\bf w}_s^{(i_1)}-{\bf w}_t^{(i_1)}\right)^{(q)}
d{\bf w}_s^{(i_2)}=
$$

\vspace{2mm}
\begin{equation}
\label{l10}
=
\frac{T-t}{2}\left(\zeta_0^{(i_1)}\zeta_0^{(i_2)}+\sum_{j=1}^{\infty}
\frac{1}{\sqrt{4j^2-1}}\left(
\zeta_{j-1}^{(i_1)}\zeta_{j}^{(i_2)}-
\zeta_j^{(i_1)}\zeta_{j-1}^{(i_2)}\right)\right).
\end{equation}

\vspace{5mm}

From (\ref{l10}) it follows that the equality (\ref{l100}) is fulfilled.
It is not difficult to see that the relation (\ref{l200}) 
can also be obtained using the approach from this section.

Let $\{\phi_j(\tau)\}_{j=0}^{\infty}$ be the complete 
orthonormal system of trigonomertic functions in the space $L_2([t, T]),$
which has the form (\ref{rre}).

We have

\vspace{-2mm}
\begin{equation}
\label{rre11}
\int\limits_t^s
\phi_j(\tau)d\tau=
\frac{T-t}{2\pi r}\left\{
\begin{matrix}
\phi_{2r-1}(s),\ &\ j=2r\cr\cr\cr
\sqrt{2}\phi_0(s)-\phi_{2r}(s),\ &\ j=2r-1
\end{matrix}
\right.\ ,
\end{equation}

\vspace{4mm}
\noindent
where $j\ge 1$ and $r=1, 2,\ldots.$

From (\ref{l1}) and (\ref{rre11}) we obtain

$$
\int\limits_t^T
\left({\bf w}_s^{(i_1)}-{\bf w}_t^{(i_1)}\right)^{(q)}
d{\bf w}_s^{(i_2)}=
\frac{1}{\sqrt{T-t}}\zeta_0^{(i_1)}
\int\limits_t^T(s-t){\bf w}_s^{(i_2)}+
$$

\vspace{1mm}
$$
+
\frac{T-t}{2}\sum_{r=1}^q
\frac{1}{\pi r}
\left(\left(\zeta_{2r}^{(i_1)}\zeta_{2r-1}^{(i_2)}-
\zeta_{2r-1}^{(i_1)}\zeta_{2r}^{(i_2)}\right)+
\sqrt{2}\zeta_{0}^{(i_2)}\zeta_{2r-1}^{(i_1)}\right)=
$$

\vspace{4mm}

$$
=
\frac{1}{\sqrt{T-t}}\zeta_0^{(i_1)}
\frac{(T-t)^{3/2}}{2}\left(\zeta_0^{(i_2)}-
\frac{\sqrt{2}}{\pi}\sum_{r=1}^{\infty}\frac{1}{r}
\zeta_{2r-1}^{(i_2)}\right)+
$$

\vspace{1mm}
$$
+
\frac{T-t}{2}\sum_{r=1}^q
\frac{1}{\pi r}
\left(\left(\zeta_{2r}^{(i_1)}\zeta_{2r-1}^{(i_2)}-
\zeta_{2r-1}^{(i_1)}\zeta_{2r}^{(i_2)}\right)+
\sqrt{2}\zeta_{0}^{(i_2)}\zeta_{2r-1}^{(i_1)}\right)=
$$

\vspace{4mm}

$$
=\frac{1}{2}(T-t)\Biggl(
\zeta_{0}^{(i_1)}\zeta_{0}^{(i_2)}
+\frac{1}{\pi}
\sum_{r=1}^{q}\frac{1}{r}\biggl(
\zeta_{2r}^{(i_1)}\zeta_{2r-1}^{(i_2)}-
\zeta_{2r-1}^{(i_1)}\zeta_{2r}^{(i_2)}+\biggr.\Biggr.
$$

\vspace{1mm}
$$
+\biggl.\Biggl.
\sqrt{2}\left(\zeta_{2r-1}^{(i_1)}\zeta_{0}^{(i_2)}-
\zeta_{0}^{(i_1)}\zeta_{2r-1}^{(i_2)}\right)\biggr)\Biggr)-
$$

\vspace{1mm}
\begin{equation}
\label{zq1}
-\frac{T-t}{\pi \sqrt{2}}
\zeta_{0}^{(i_1)}\sum_{r=q+1}^{\infty}\frac{1}{r}
\zeta_{2r-1}^{(i_2)}.
\end{equation}

\vspace{6mm}

From (\ref{zq1}) and (\ref{rez1}) we obviously get
(\ref{ajjja}).

\vspace{5mm}

\section{Convergence in the 
Mean of Degree $2n$ and With Probability 1}

\vspace{5mm}

Let us denote

\vspace{-1mm}
\begin{equation}
\label{555}
A_{T,t}^{(i_1 i_2)q}=
\frac{T-t}{2}\sum_{i=1}^{q}
\frac{1}{\sqrt{4i^2-1}}\left(
\zeta_{i-1}^{(i_1)}\zeta_{i}^{(i_2)}-
\zeta_i^{(i_1)}\zeta_{i-1}^{(i_2)}\right),
\end{equation}

\vspace{1mm}
\begin{equation}
\label{abc}
{\hat A}_{T,t}^{(i_1 i_2)q}=
\frac{T-t}{2\pi}
\sum_{r=1}^{q}\frac{1}{r}\left(
\zeta_{2r}^{(i_1)}\zeta_{2r-1}^{(i_2)}-
\zeta_{2r-1}^{(i_1)}\zeta_{2r}^{(i_2)}+
\sqrt{2}\left(\zeta_{2r-1}^{(i_1)}\zeta_{0}^{(i_2)}-
\zeta_{0}^{(i_1)}\zeta_{2r-1}^{(i_2)}\right)\right).
\end{equation}

\vspace{6mm}

Then, from (\ref{2121}) we get

\vspace{-1mm}
\begin{equation}
\label{bb1}
I_{T,t}^{(i_1 i_2)q}=\frac{T-t}{2}\zeta_{0}^{(i_1)}\zeta_{0}^{(i_2)}+
A_{T,t}^{(i_1 i_2)q},
\end{equation}

\begin{equation}
\label{bb2}
I_{T,t}^{(i_1 i_2)q}=
\frac{T-t}{2}\zeta_{0}^{(i_1)}\zeta_{0}^{(i_2)}+
{\hat A}_{T,t}^{(i_1 i_2)q}.
\end{equation}

\vspace{4mm}

It is not difficult to demonstrate \cite{Mi2} that from (\ref{ajjja})
we can get an another representation for the Levy stochastic area

$$
\frac{T-t}{2\pi}\Biggl(
\sum_{r=1}^{q}\frac{1}{r}\left(
\zeta_{2r}^{(i_1)}\zeta_{2r-1}^{(i_2)}-
\zeta_{2r-1}^{(i_1)}\zeta_{2r}^{(i_2)}+
\right.\Biggr.
$$

\vspace{2mm}
$$
+\Biggl.\left.\sqrt{2}\left(\zeta_{2r-1}^{(i_1)}\zeta_{0}^{(i_2)}-
\zeta_{0}^{(i_1)}\zeta_{2r-1}^{(i_2)}\right)\right)
+\sqrt{2}\left(\frac{\pi^2}{6}-\sum_{r=1}^q\frac{1}{r^2}\right)^{1/2}\left(
\xi_q^{(i_1)}\zeta_0^{(i_2)}-\zeta_0^{(i_1)}\xi_q^{(i_2)}\right)\Biggr),
$$

\vspace{5mm}
\noindent
where

\vspace{-3mm}
$$
\xi_q^{(i)}=\left(\frac{\pi^2}{6}-
\sum\limits_{r=1}^q\frac{1}{r^2}\right)^{-1/2}
\sum_{r=q+1}^{\infty}
\frac{1}{r}\zeta_{2r-1}^{(i)},
$$

\vspace{5mm}
\noindent
and $\zeta_0^{(i)},$ $\zeta_{2r}^{(i)},$
$\zeta_{2r-1}^{(i)},$ $\xi_q^{(i)}$ ($r=1,\ldots,q,$\
$i=1,\ldots,m$) are independent standard Gaussian random variables.

From (\ref{555}) and (\ref{abc}) we obtain

\vspace{2mm}
$$
{\sf M}\left\{\left(A_{T,t}^{(i_1 i_2)}- A_{T,t}^{(i_1 i_2)q}
\right)^2\right\}=
\frac{(T-t)^2}{2}\Biggl(\frac{1}{2}-\sum_{i=1}^q
\frac{1}{4i^2-1}\Biggr)=
$$

\vspace{2mm}
$$
=
\frac{(T-t)^2}{2}
\sum\limits_{i=q+1}^{\infty}\frac{1}{4i^2-1}
\le \frac{(T-t)^2}{2}\int\limits_{q}^{\infty}
\frac{1}{4x^2-1}dx
=
$$

\vspace{2mm}
\begin{equation}
\label{teac}
=-\frac{(T-t)^2}{8}{\rm ln}\left|
1-\frac{2}{2q+1}\right|
\le C_1\frac{(T-t)^2}{q},
\end{equation}

\vspace{6mm}

$$
{\sf M}\left\{\left({\hat A}_{T,t}^{(i_1 i_2)}-{\hat A}_{T,t}^{(i_1 i_2)q}
\right)^2\right\}
=
\frac{3(T-t)^{2}}{2\pi^2}
\left(\frac{\pi^2}{6}-\sum\limits_{r=1}^q\frac{1}{r^2}\right)
=
$$

\vspace{2mm}
$$
=
\frac{3(T-t)^{2}}{2\pi^2}\sum\limits_{r=q+1}^{\infty}\frac{1}{r^2}
\le\frac{3(T-t)^{2}}{2\pi^2}\int\limits_q^{\infty}
\frac{dx}{x^2}
=
$$

\vspace{2mm}

\begin{equation}
\label{teac0}
=\frac{3(T-t)^{2}}{2\pi^2 q}
\le C_2 \frac{(T-t)^2}{q},
\end{equation}

\vspace{6mm}
\noindent
where constants $C_1,$ $C_2$ does not depend on $q.$
\vspace{2mm}

For the case $k=2,$ $i_1\ne i_2$, and $\psi_1(s),$ $\psi_2(s)\equiv 1$ 
from (\ref{99999}) we obtain

\vspace{1mm}
\begin{equation}
\label{yyy}
{\sf M}\left\{\left(I_{T,t}^{(i_1 i_2)}-I_{T,t}^{(i_1 i_2)q}
\right)^{2n}\right\}\le C_{n,2}\left(\frac{(T-t)^2}{2}
\left(\frac{1}{2}-\sum_{i=1}^q
\frac{1}{4i^2-1}\right)\right)^n\ \to 0\ \ \ \hbox{if}\ q\to\infty,
\end{equation}

\vspace{2mm}

\begin{equation}
\label{yyy2}
{\sf M}\left\{\left(I_{T,t}^{(i_1 i_2)}-I_{T,t}^{(i_1 i_2)q}
\right)^{2n}\right\}\le C_{n,2}
\left(\frac{3(T-t)^{2}}{2\pi^2}\left(\frac{\pi^2}{6}-
\sum\limits_{r=1}^q\frac{1}{r^2}\right)
\right)
^n\ \to 0\ \ \ \hbox{if}\ q\to\infty,
\end{equation}

\vspace{5mm}
\noindent
where
$C_{n,k}=(k!)^{n} (2n-1)^{nk},$
$I_{T,t}^{(i_1 i_2)q}$ 
has the form  (\ref{bb1}) in the inequality (\ref{yyy}), 
and $I_{T,t}^{(i_1 i_2)q}$
has the form (\ref{bb2}) in the inequality (\ref{yyy2}), 

From (\ref{teac})--(\ref{yyy2}) we get

\vspace{1mm}
$$
{\sf M}\left\{\left(A_{T,t}^{(i_1 i_2)}-A_{T,t}^{(i_1 i_2)q}
\right)^{2n}\right\}\ \to 0\ \ \ {\rm if}\ \ \ q\to\infty,
$$

\vspace{1mm}
$$
{\sf M}\left\{\left(\hat A_{T,t}^{(i_1 i_2)}-\hat A_{T,t}^{(i_1 i_2)q}
\right)^{2n}\right\}\ \to 0\ \ \ {\rm if}\ \ \ q\to\infty.
$$

\vspace{4mm}

Let us address now to the convergence w.~p.~1
for $A_{T,t}^{(i_1 i_2)q}$. 
First, note the well known fact.

\vspace{2mm}

{\bf Lemma 1.}\ {\it If for the sequence of random variables
$\xi_q$ and for some
$\alpha>0$ the number series 

$$
\sum\limits_{q=1}^{\infty}{\sf M}\left\{|\xi_q|^{\alpha}\right\}
$$

\vspace{3mm}
\noindent
converges, then the sequence $\xi_q$ converges to zero w.~p.~{\rm 1}.}

\vspace{2mm}

From (\ref{teac}) and (\ref{yyy}) $(n=2)$ we obtain

\vspace{1mm}
$$
{\sf M}\left\{\left(I_{T,t}^{(i_1 i_2)}-I_{T,t}^{(i_1 i_2)q}
\right)^{4}\right\}={\sf M}
\left\{\left(A_{T,t}^{(i_1 i_2)}-A_{T,t}^{(i_1 i_2)q}
\right)^{4}\right\}\le \frac{K}{q^2},
$$

\vspace{4mm}
\noindent
where constant $K$ does not depend on $q.$

Since the series

\vspace{-1mm}
$$
\sum\limits_{q=1}^{\infty}\frac{K}{q^2}
$$

\vspace{3mm}
\noindent
converges, then according to Lemma 1
we obtain that
$A_{T,t}^{(i_1 i_2)}-A_{T,t}^{(i_1 i_2)q}\to 0$\ if\  $q\to \infty$ 
w.~p.~1. Then  
$A_{T,t}^{(i_1 i_2)q}\to A_{T,t}^{(i_1 i_2)}$\  if\  $q\to \infty$ 
w.~p.~1.

In addition, using (\ref{teac0}) and (\ref{yyy2}) $(n=2)$, we get
$\hat A_{T,t}^{(i_1 i_2)q}\to \hat A_{T,t}^{(i_1 i_2)}$\  if\  $q\to \infty$ 
w.~p.~1.

\end{document}